\documentclass[12pt]{article}
\usepackage{graphicx}
\usepackage{amsmath}
\usepackage{amsfonts}
\usepackage{amssymb}
\newtheorem{theorem}{Theorem}

\newtheorem{lemma}[theorem]{Lemma}

\newtheorem{proposition}[theorem]{Proposition}

\newenvironment{proof}[1][Proof]{\textbf{#1.} }{\ \rule{0.5em}{0.5em}}

\def\text{\hbox} \def\newpage{\vfill\break}

\def\a{\alpha}
\def\b{\beta}

\def\p{\pi}

\def\r{\rho}

\def\s{\sigma}

\def\L{\Lambda}

\def\S{\Sigma}

\def\GG{{\mathcal G}}

\def\L{{\mathcal L}}
\def\Z{{\mathbf Z}}
\def\S{{\mathbf S}}
\def\Q{{\mathbf Q}}
\def\QQ{\widetilde{{\mathbf Q}}}

\begin{document}

\title{On periodic Takahashi manifolds\footnote{Work performed under the auspices of
G.N.S.A.G.A. of C.N.R. of Italy and supported by the University of
Bologna, funds for selected research topics.}}
\author{Michele Mulazzani}

\maketitle
\begin{abstract}
{In this paper we show that periodic Takahashi 3-manifolds are
cyclic coverings of the connected sum of two lens spaces (possibly
cyclic coverings of $\S^3$), branched over knots. When the base
space is a 3-sphere, we prove that the associated branching set is
a two-bridge knot of genus one, and we determine its type.
Moreover, a geometric cyclic presentation for the fundamental
groups of these manifolds is obtained in several interesting
cases, including the ones corresponding to the branched cyclic
coverings of $\S^3$.
}

\end{abstract}


\section{Introduction}
\begin{figure}[ht]
 \begin{center}
 \includegraphics*[totalheight=3cm]{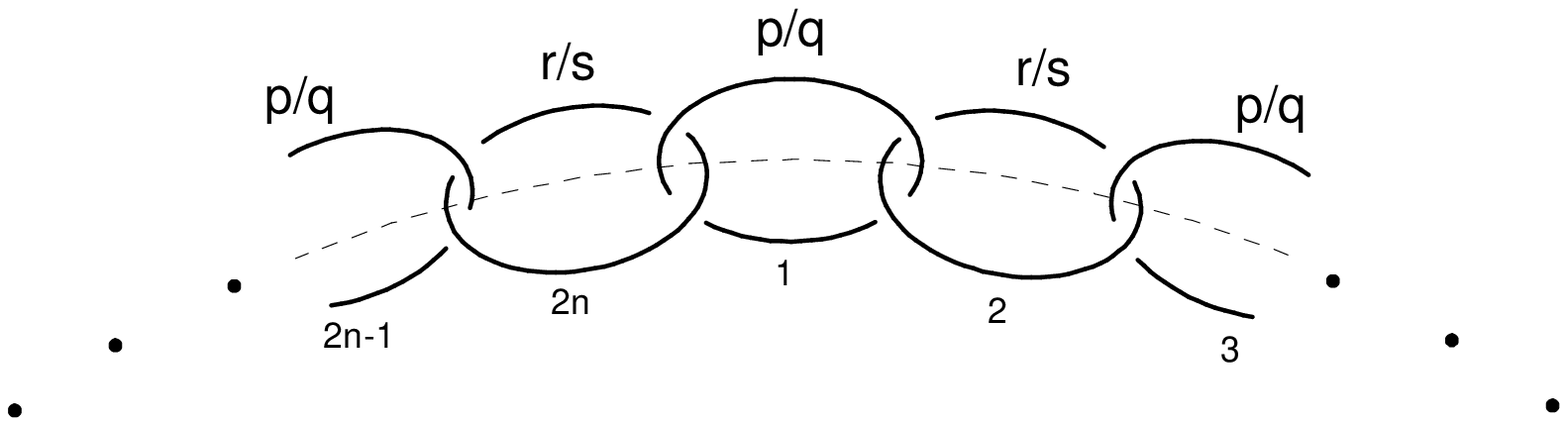}
 \end{center}
 \caption{Surgery along ${\cal L}_{2n}$ yielding $M_n(p/q,r/s)$.}

 \label{Fig. 1}

\end{figure}

Takahashi manifolds are closed orientable 3-manifolds introduced
in \cite{Ta} by Dehn surgery on $\S^3$, with rational
coefficients, along the $2n$-component link ${\cal L}_{2n}$
depicted in Figure \ref{Fig. 1}. These manifolds have been
intensively studied in \cite{KKV}, \cite{RS}, and \cite{VK}. In
the latter two papers, a nice topological characterization of all
Takahashi manifolds as two-fold coverings of $\S^3$, branched over
the closure of certain rational 3-string braids, is given.

A Takahashi manifold is called {\it periodic\/} when the surgery
coefficients have the same cyclic symmetry of order $n$ of the
link ${\cal L}_{2n}$, i.e. the coefficients are $p/q$ and $r/s$
alternately. Several important classes of 3-manifolds, such as
(fractional) Fibonacci manifolds \cite{HKM,VK} and
Sieradsky manifolds \cite{CHK,Si}, represent notable
examples of periodic Takahashi manifolds.

In this paper we show that each periodic Takahashi manifold is an
$n$-fold cyclic covering of the connected sum of two lens spaces,
branched over a knot. This knot arises from a component of the Borromean
rings, by performing a surgery with coefficients $p/q$ and $r/s$ along the
other two components.

For particular values of the surgery coefficients (including the
classes of manifolds cited above), the periodic Takahashi
manifolds turn out to be $n$-fold cyclic coverings of $\S^3$,
branched over two-bridge knots of genus one\footnote{For notation and
properties about two-bridge knots and links we refer to
\cite{BZ}. For the characterization of
two-bridge knots of genus one, see \cite{Fu}.}, whose parameters
are obtained using Kirby-Rolfsen calculus \cite{Ro} (compare the analogous result of
\cite{KKV}, obtained by a different approach). Observe that in
\cite {RS} a characterization of all periodic Takahashi manifolds
as $n$-fold cyclic coverings of $\S^3$, branched over the closure
of certain rational 3-string braids, is presented, but the result
is incorrect, as we show in Remark 1.

For many interesting periodic Takahashi manifolds - including the
ones corresponding to branched cyclic coverings of $\S^3$ - a
cyclic presentation for the fundamental group is provided and
proved to be geometric, i.e. arising from a Heegaard diagram, or,
equivalently, from a canonical spine\footnote{A canonical spine is
a 2-dimensional cell complex with a single vertex.} \cite{Ne}.

\section{Main results}

We denote by $M(p_1/q_1,\ldots,p_n/q_n;r_1/s_1,\ldots,r_n/s_n)$
the Takahashi manifold obtained by Dehn surgery on $\S^3$ along
the $2n$-component link ${\cal L}_{2n}$ of Figure \ref{Fig. 1},
with surgery coefficients
$p_1/q_1,r_1/s_1,\ldots,p_n/q_n,r_n/s_n\in\QQ= \Q\cup\{\infty\}$
respectively, cyclically associated to the components of ${\cal
L}_{2n}$.

A Takahashi manifold is periodic when $p_i/q_i=p/q$ and
$r_i/s_i=r/s$, for every $i=1,\ldots,n$. Denote by $M_n(p/q,r/s)$
the periodic Takahashi manifold
$M(p/q,\ldots,p/q;r/s,\ldots,r/s)$. From now on, without loss of
generality, we can always suppose that: $\gcd(p,q)=1$,
$\gcd(r,s)=1$ and $p,r\ge 0$. Moreover, if $\a,\b\in{\bf Z}$ with
$\a\ge 0$ and $\gcd(\a,\b)=1$, we shall denote by $L(\a,\b)$ the
lens space of type $(\a,\b)$. As usual, $L(0,1)$ is homeomorphic
to $\S^1\times\S^2$ and $L(1,\b)$ is homeomorphic to $\S^3$, for
all $\b$ (including $\b=0$).

Notice that $M_n(p/q,-p/q)$ is the Fractional Fibonacci manifold
$M_n^{p/q}$ defined in \cite{VK} and, in particular, $M_n(1,-1)$
is the Fibonacci manifold $M_n$ studied in \cite{HKM}. Moreover,
$M_n(1,1)$ is the Sieradsky manifold $M_n$ introduced in \cite{Si}
and studied in \cite{CHK}. Because of the symmetries of ${\cal
L}_{2n}$, the homeomorphisms $$M_n(p/q,r/s)\cong
M_n(-p/q,-r/s)\cong M_n(r/s,p/q)\cong M_n(-r/s,-p/q)$$ can easily
be obtained for all $n\ge 1$ and $p/q,r/s\in\QQ$.

A balanced presentation of the fundamental group of every
Takahashi manifold is given in \cite{Ta}, and in \cite{RS} it is
shown that this presentation is geometric, i.e. it arises from a
Heegaard diagram (or, equivalently, from a canonical spine). As a
consequence, $\p_1(M_n(p/q,r/s))$ admits the following geometric
presentation with $2n$ generators and $2n$ relators:
$$<x_1,\ldots,x_{2n}\,\vert\,
x_{2i-1}^qx_{2i}^{-r}x_{2i+1}^{-q},\,
x_{2i}^sx_{2i+1}^px_{2i+2}^{-s}\,;\,i=1,\ldots,n>,$$ where the
subscripts are mod $2n$.

When $r=1$, we can easily get a cyclic presentation  \cite{Jo}
with $n$ generators:\footnote{Alternatively, a similar cyclic
presentation can be obtained when $p=1$.} $$\p_1(M_n(p/q,1/s))=
<z_1,\ldots,z_n\,\vert\,z_i^p(z_i^{-q}z_{i+1}^q)^s(z_i^{-q}z_{i-1}^q)^s\,
;\,i=1,\ldots,n>,\eqno(1)$$ where the subscripts are mod $n$.

\begin{proposition} \label{Proposition new}
For all $p/q\in\QQ$ and $s\in\bf Z$, the cyclic presentation (1)
of $\p_1(M_n(p/q,1/s))$ is geometric.
\end{proposition}

\begin{figure}
 \begin{center}
 \includegraphics*[totalheight=18cm]{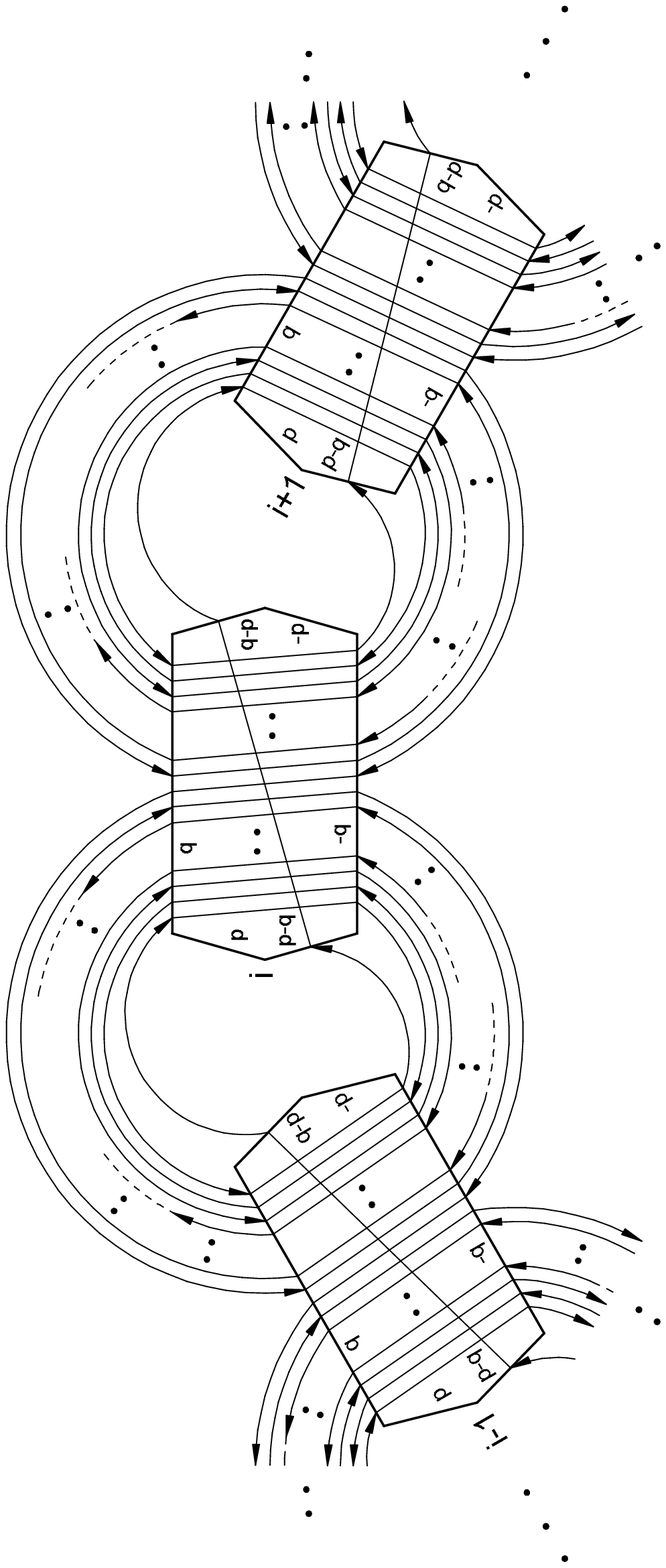}
 \end{center}
 \caption{An RR-system for the cyclic presentation ($1'$).}
 \label{Fig. 2}
\end{figure}

\begin{proof}
If $s=0$ then $M_n(p/q,1/s)$ is homeomorphic to the connected sum
of $n$ copies of $L(p/q)$, and therefore the statement is
straightforward. If $s>0$, the presentation becomes
$$<z_1,\ldots,z_n\,\vert\,z_i^{p-q}(z_{i+1}^qz_i^{-q})^s(z_{i-1}^qz_i^{-q})^{s-1}z_{i-1}^q\,
;\,i=1,\ldots,n>.\eqno(1')$$ Figure \ref{Fig. 2} shows an
RR-system which induces ($1'$), and so, by \cite{OS}, this
presentation is geometric. If $s<0$, the presentation becomes
$$<z_1,\ldots,z_n\,\vert\,z_i^{p+q}(z_{i+1}^{-q}z_i^q)^{-s}(z_{i-1}^{-q}z_i^q)^{-s-1}z_{i-1}^{-q}\,
;\,i=1,\ldots,n>.\eqno(1'')$$ Therefore,  if we replace $q$ with
$-q$, Figure \ref{Fig. 2} also gives an RR-system inducing
($1''$).
\end{proof}

\medskip
Since the link ${\cal L}_2$ is a two-component trivial link, we
immediately get the following results:

\begin{lemma} \label{Lemma 1} For all $p/q,r/s\in\QQ$,
the manifold $M_1(p/q,r/s)$ is homeomorphic to the connected sum
of lens spaces $L(p,q)\#L(r,s)$. In particular, $M_1(p/q,1/s)$ is
homeomorphic to the lens space $L(p,q)$ and $M_1(1/q,1/s)$  is
homeomorphic to $\S^3$.
\end{lemma}

\begin{proof}
$M_1(p/q,r/s)$ is obtained by Dehn surgery on $\S^3$, with
coefficients $p/q$ and $r/s$, along the trivial link with two
components ${\cal L}_{2}$.
\end{proof}

\medskip

Now we prove the main result of the paper:

\begin{theorem} \label{main theorem} For all $p/q,r/s\in\QQ$
and $n>1$, the periodic Takahashi manifold $M_n(p/q,r/s)$ is the
$n$-fold cyclic covering of the connected sum of lens spaces
$L(p,q)\#L(r,s)$, branched over a knot $K$ which does not depend
on $n$. Moreover, $K$ arises from a component of the Borromean
rings, by performing a surgery with coefficients $p/q$ and $r/s$ along the
other two components.
\end{theorem}

\newpage

\begin{figure}
 \begin{center}
 \includegraphics*[totalheight=6cm]{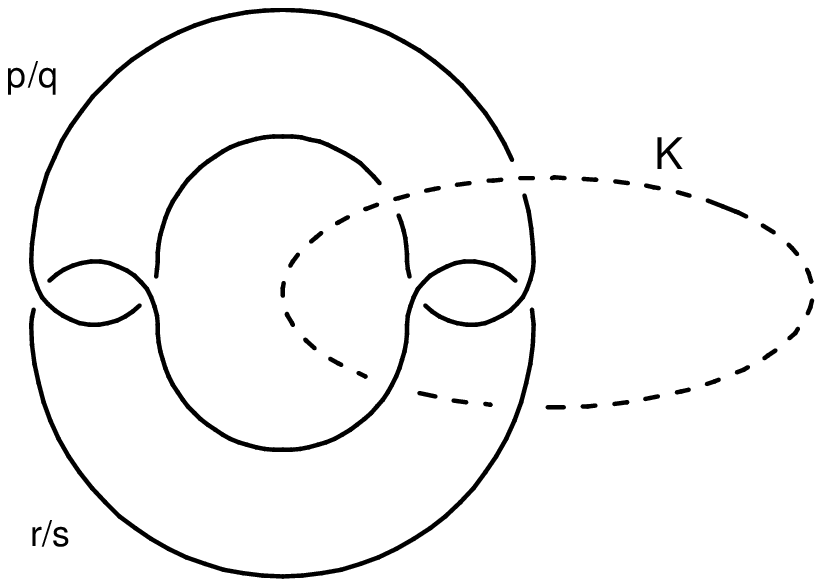}
 \end{center}
 \caption {The branching set $K$ (dashed line).}

 \label{Fig. 0}

\end{figure}

\begin{proof}
Both the link ${\cal L}_{2n}$ and the surgery coefficients
defining $M_n(p/q,r/s)$ are invariant with respect to the rotation
$\r_n$ of $\S^3$, which sends the $i$-th component of ${\cal
L}_{2n}$ onto the $(i+2)$-th component (mod $2n$). Let ${\cal
G}_n$ be the cyclic group of order $n$ generated by $\r_n$.
Observe that the fixed-point set of the action of ${\cal G}_n$ on
$\S^3$ is a trivial knot disjoint from ${\cal L}_{2n}$. Therefore,
we have an action of ${\cal G}_n$ on $M_n(p/q,r/s)$, with a knot
$K_n$ as fixed-point set. The quotient $M_n(p/q,r/s)/{\cal G}_n$
is precisely the manifold $M_1(p/q,r/s)$, which is homeomorphic to
$L(p,q)\#L(r,s)$ by Lemma \ref{Lemma 1}, and $K_n/{\cal G}_n$ is
obviously a knot $K\subset M_1(p/q,r/s)$, which only depends on
$p/q$ and $r/s$. Moreover, $K\cup{\cal L}_{2}$ is the Borromean
rings, as showed in Figure \ref{Fig. 0}. This proves the statement.
\end{proof}

\medskip

We can give another description of the branching set $K$, as
the inverse image of a trivial knot in a certain two-fold branched
covering.

\begin{figure}[ht]
 \begin{center}
 \includegraphics*[totalheight=5cm]{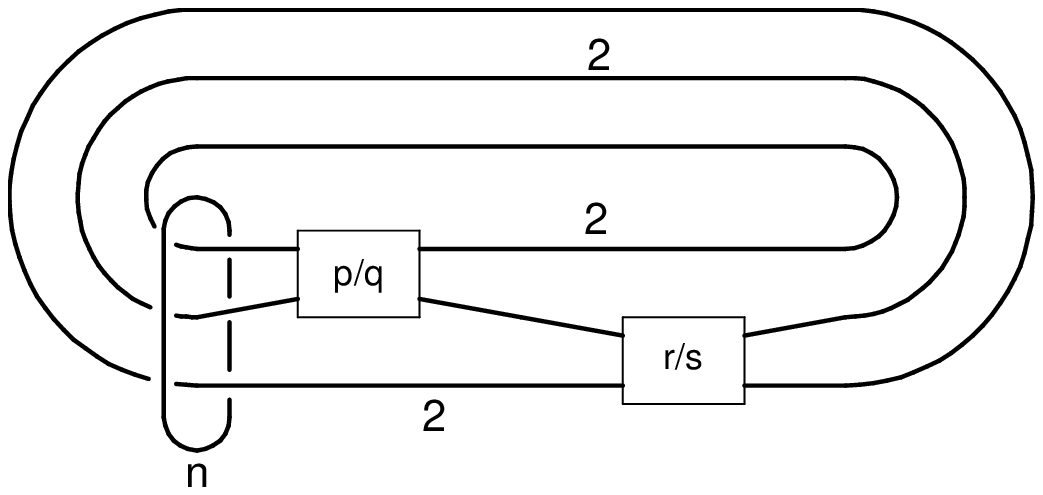}
 \end{center}
 \caption {The link $\L(p/q,r/s)$.}

 \label{Fig. 00}

\end{figure}

Denote by $\L(p/q,r/s)$ the link depicted in Figure \ref{Fig. 00}.
It is composed by the closure of the rational 3-string braid
$\s_{1}^{p/q}\s_2^{r/s}$, which is the connected sum of the
two-bridge knots or links ${\bf b}(p,q)$ and ${\bf b}(r,s)$, and
by a trivial knot. Moreover, denote: (i) by ${\cal O}_{n} (p / q ,
r / s )= M_n (p/ q , r / s ) / {\cal G}_n$ the orbifold from the
proof of Theorem \ref{main theorem}, whose underlying space is
$L(p,q)\#L(r,s)$ and whose singular set is the knot $K$, with
index $n$; (ii) by $\S^3({\cal K}_{n} (p / q , r / s ))$ the
orbifold whose underlying space is $\S^3$ and whose singular set
is the closure of the rational 3-string braid
$(\s_{1}^{p/q}\s_2^{r/s})^n$, with index $2$; and (iii) by
$\S^3(\L(p/q,r/s))$ the orbifold  whose underlying space is $\S^3$
and whose singular set is the link $\L(p/q,r/s)$, with
index $2$ and $n$ as pointed out in Figure \ref{Fig.
00}.

\begin{proposition} \label{commutative}
Assuming the previous notations, the following commutative diagram
holds for each periodic Takahashi manifold.

 \begin{center}
  \unitlength=0.5mm
  \begin{picture}(90,75)(0,10)
  \put(40,80){\makebox(0,0)[cc]{$M_{n} (p/ q , r / s )$}}
  \put(50,70){\vector(2,-1){30}}
  \put(30,70){\vector(-2,-1){30}}
  \put(68,68){\makebox(0,0)[cc]{$n$}} 
  \put(7,68){\makebox(0,0)[cc]{$2$}} 
  \put(0,45){\makebox(0,0)[cc]{$\S^3({\cal K}_{n} (p/ q , r/ s ))$}}
  \put(80,45){\makebox(0,0)[cc]{${\cal O}_{n} (p / q , r / s )$}}
  \put(0,35){\vector(2,-1){30}}
  \put(80,35){\vector(-2,-1){30}}
  \put(7,25){\makebox(0,0)[cc]{$n$}} 
  \put(68,25){\makebox(0,0)[cc]{$2$}} 
  \put(40,10){\makebox(0,0)[cc]{$\S^3(\L(p/q,r/s))$}}
  \end{picture}
  \end{center}

\end{proposition}

\begin{proof} The link ${\cal L}_{2n}$ admits an
invertible involution $\tau$, whose axis intersects each component
in two points (see the dashed line of Figure \ref{Fig. 1}), and
the rotation symmetry $\r_n$ of order $n$ which was discussed in
Theorem \ref{main theorem}. These symmetries induce symmetries
(also denoted by $\tau$ and $\r_n$) on the periodic Takahashi
manifold $M = M_n(p/ q , r / s )$, such that $\langle \tau,\r_n
\rangle\cong\langle \tau \rangle \oplus \GG_n\cong \Z_2 \oplus
\Z_n$. We have $M/\langle \tau \rangle=\S^3({\cal K}_{n} (p/ q ,
r/ s ))$ (see \cite{RS} and \cite{VK}) and $M/\GG_n={\cal O}_{n}
(p / q , r / s )$ (see Theorem \ref{main theorem}). It is
immediate to see that $\r_n$ induces a symmetry (also denoted by
$\r_n$) on the orbifold $M/ \langle \tau \rangle $, and $(M/
\langle \tau \rangle )/\GG_n$ is the orbifold $\S^3(\L(p/q,r/s))$.
As we see from Figure \ref{Fig. 0}, $\tau$ induces a strongly
invertible involution (also denoted by $\tau$) on the link ${\cal
L}_{2}$. Using the Montesinos algorithm we see that $(M / \GG_n) /
\langle \tau \rangle = \S^3(\L(p/q,r/s))$. This concludes the
proof.
\end{proof}

\medskip

As a consequence, the branching set $K$ of Theorem \ref{main
theorem} can be obtained as the inverse image of the trivial
component of $\L(p/q,r/s)$ in the two-fold branched covering
${\cal O}_{n} (p / q , r / s )\to\S^3(\L(p/q,r/s))$.

\medskip

From Theorem \ref{main theorem} we can get the following result,
which has already been obtained in \cite{KKV} by a different
technique.

\begin{proposition} \label{Proposition 1} For all $q,s\in\bf Z$ and $n>1$,
the periodic Takahashi manifold $M_n(1/q,1/s)$ is the $n$-fold
cyclic covering of $\S^3$, branched over the two-bridge knot of
genus one ${\bf b}(\vert 4sq-1\vert,2s)\cong{\bf b}(\vert
4sq-1\vert,2q)$.
\end{proposition}

\begin{figure}
 \begin{center}
 \includegraphics*[totalheight=14cm]{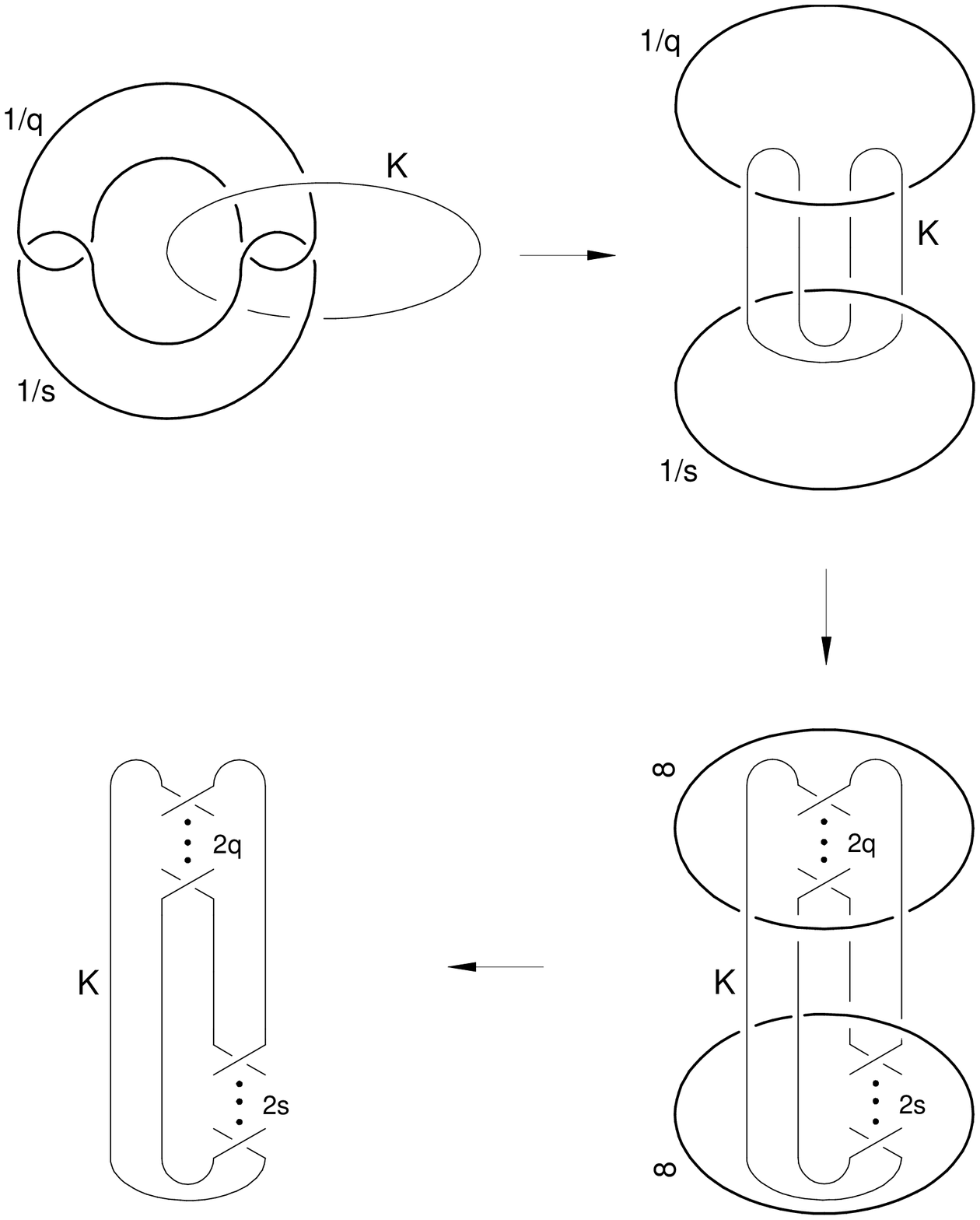}
 \end{center}
 \caption{}
 \label{Fig. 10}
\end{figure}

\begin{proof}
From Theorem \ref{main theorem}, $M_n(1/q,1/s)$ is the $n$-fold
cyclic covering of $L(1,q)\#L(1,s)\cong\S^3$, branched over a knot
$K$ which does not depend on $n$. By isotopy and Kirby-Rolfsen moves it is
easy to obtain (see Figure \ref{Fig. 10}) a diagram of $K$, which
is a Conway's normal form of type $[-2q,2s]$. This proves the
statement.
\end{proof}

\medskip

Proposition \ref{Proposition 1} covers the results of \cite{CHK},
\cite{HKM} and \cite{VK} concerning $n$-fold branched cyclic
coverings of two-bridge knots. Moreover, for all $p,q\in\bf Z$,
the periodic Takahashi manifold $M_n(1/q,1/s)$ is homeomorphic to
the Lins-Mandel manifold $S(n,\vert 4sq-1\vert,2s,1)$
\cite{LM,Mu}, the Minkus manifold $M_n(\vert 4sq-1\vert,2s)$
\cite{Mi} and the Dunwoody manifold $M((\vert
4q-1\vert-1)/2,0,1,s,n,-q_{\s})$ \cite{Du,GM}.

Moreover, observe that all cyclic coverings of two-bridge knots of
genus one are periodic Takahashi manifolds.

\medskip

\noindent{\bf Remark 1} The results of Corollaries 8, 9 and 11 of
\cite{RS}, concerning periodic Takahashi manifolds as $n$-fold
cyclic branched coverings of the closure of certain (rational)
3-string braids, are incorrect. This is evident from the following
counterexamples. If $p/q=3$ and $r/s=-3$ then the first homology
group of the 3-fold cyclic branched covering of the closure of the
3-string braid $(\s_1^3\s_2^{-3})^2$ has order $256$, but $\vert
H_1(M_3(3,-3))\vert=1296$. If $p/q=3/2$ and $r/s=1$ then the first
homology group of the 4-fold cyclic branched covering of the
closure of the rational 3-string braid $(\s_1^{3/2}\s_2)^2$ has
order $135$, but $\vert H_1(M_4(3/2,1))\vert=15$. Note that the
corollaries are valid if $p=r=1$.

\medskip

The following conjecture is naturally suggested by the previous results.

\medskip

\noindent{\bf Conjecture} Let $p/q,r/s\in\QQ$ be fixed. Then, for
all $n>1$, the periodic Takahashi manifolds $T_n=M_n(p/q,r/s)$ are
$n$-fold cyclic coverings of $\S^3$, branched over a knot which
does not depend on $n$, if and only if $p=1=r$.

\bigskip

\noindent \underline{Added in revision} - The referee pointed out
that it is possible to prove the conjecture for ``almost all
cases'' by using the hyperbolic Dehn surgery theorem and the shortest
geodesic arguments by Kojima \cite{Ko}.

\bigskip

\noindent{\bf Acknowledgement}

\noindent The author wishes to thank the referee for his valuable
suggestions to improve this paper and Massimo Ferri and Andrei
Vesnin for the useful discussions on the topics.


\vspace{15 pt} {MICHELE MULAZZANI, Department of Mathematics,
University of Bologna, I-40127 Bologna, ITALY, and C.I.R.A.M.,
Bologna, ITALY. E-mail: mulazza@dm.unibo.it}

\end{document}